\documentclass[11pt]{article}

\usepackage[margin = 1.0in]{geometry}   
\usepackage{geometry}  			
\usepackage{amssymb}
\usepackage{amsmath}
\usepackage{palatino}
\usepackage{graphicx}
\usepackage[english]{babel}
\usepackage{amscd}
\usepackage{amsthm}
\usepackage{color}
\usepackage{enumerate}
\usepackage{hyperref}
\usepackage{subcaption}
\usepackage{comment}
\newtheorem{theorem}{Theorem}
\newtheorem{conjecture}{Conjecture}
\newtheorem{corollary}[theorem]{Corollary}
\newtheorem{lemma}{Lemma}

\newcommand{\Prb}{\mathbb{P}}
\newcommand{\Exp}{\mathbb{E}}
\newcommand{\mbI}{\mathbb{I}}

\newcommand{\bh}{\mathbf{H}}

\title{An entropic analogue of the MMS conjecture}

\author{
Jianhang Ai\thanks{Faculty of Electrical Engineering, Czech Technical University, Karlovo N\'am\v{e}st\'i 13, 12135, Prague 2, Czech Republic, e-mail: ai.jianhang@fel.cvut.cz}
\and
Ond\v{r}ej Ku\v{z}elka\thanks{Faculty of Electrical Engineering, Czech Technical University, Karlovo N\'am\v{e}st\'i 13, 12135, Prague 2, Czech Republic, e-mail: ondrej.kuzelka@fel.cvut.cz}
\and
Christos Pelekis\thanks{Aristotle University of Thessaloniki, 
Department of Mathematics, 
541 24 Thessaloniki, Greece, e-mail:  cpelekis@math.auth.gr }
}

\begin{document}

\maketitle

\begin{abstract}  
Let $P=\{x_1,\ldots,x_n\}$ be a multiset consisting of $n\ge 2$ real numbers such that $\sum_{i=1}^{n}x_i=0$ and $\sum_{i=1}^{n}|x_i|>0$, and let $k <n$ be a   positive integer. We sample $k$ elements from $P$ without replacement and set $X_P$ be the sum of the elements in our sample. 
It is shown  that the Shannon entropy of $X_P$ satisfies 
\[
\mathbf{H}(X_P) \ge \mathbf{H}(\text{Ber}(k/n)) \, ,
\]  
where $\text{Ber}(k/n)$ is a Bernoulli random variable of mean $k/n$.  
The result is sharp, and may be seen as an entropic analogue of the Manickam-Mikl\'os-Singhi (MMS) conjecture. 
\end{abstract}

\vspace{0.2cm}


\section{Introduction, related work and main results}

\subsection{Notation and definitions}

Let us begin with introducing some piece of notation that will be fixed throughout the text. 
Given a positive integer $n$, we denote by $[n]$ the set $\{1,\ldots,n\}$ and, given a finite set $F$, we denote by $|F|$ its cardinality. 
For $k\in \{0,1,\ldots,n\}$, we denote by $\binom{[n]}{k}$ the collection consisting of all subsets of $[n]$ having $k$ elements. 

If $X,Y$ are two random variables, we write $X\sim Y$ when they have the same distribution. 
Unless mentioned otherwise, every random variable $X$ that will be considered takes values in a \emph{finite} set $\mathcal{S}_X$. We  denote by $p_X$ the \emph{probability mass function} of the random variable $X$, that is,   $p_X(x) =\Prb(X=x), x\in\mathcal{S}_X$. 
A \emph{mode} of the random variable $X$ is a value $x^*\in\mathcal{S}_X$ such that $p_X(x^*)\ge p_X(x)$, for all $x\in\mathcal{S}_X$. 

Given positive integers 
$n,k$ such that $n\ge k\ge 1$ and a positive integer $i\in [n]$, we denote by $\text{Hyp}(n,i,k)$ a hypergeometric random variable which counts the number of black marbles in a sample without replacement of size $k$ from an urn that contains $i$
black and $n - i$ white marbles. 
Formally, if $H\sim\text{Hyp}(n,i,k)$ then  
\[
\Prb(H=j)=\frac{\binom{i}{j}\binom{n-i}{k-j}}{\binom{n}{k}}\, , \, \text{ for }  \, j\in\{\max\{0,k-(n-i)\},\dots,\min\{k,i\}\} \, . 
\]

Given a vector $\mathbf{x}= (x_1,\ldots,x_n)\in\mathbb{R}^n$, let  $x_{(1)} \ge \cdots \ge x_{(n)}$ denote the coordinates of $\mathbf{x}$ in decreasing order. 
Given two vectors $\mathbf{x}= (x_1,\ldots,x_n)$ and $\mathbf{y}= (y_1,\ldots,y_n)$  in $\mathbb{R}^n$, we say that $\mathbf{x}$ \textit{is majorized} 
by $\mathbf{y}$ if it holds 
\[
\sum_{j=1}^{\ell} x_{(j)} \le \sum_{j=1}^{\ell} y_{(j)}\, , \,  \text{ for all } \, \ell=1,\ldots,n-1 \, , \quad \text{ and } \quad \sum_{j=1}^{n} x_{(j)} = \sum_{j=1}^{n} y_{(j)} \, . 
\]
We write $\mathbf{x}\prec\mathbf{y}$ when the vector  $\mathbf{x}$  is majorized by the vector $\mathbf{y}$. 
A function $f:A\to \mathbb{R}$ defined on a set $A\subseteq\mathbb{R}^n$ is said to be \emph{Schur concave on $A$} if for every $\mathbf{x},\mathbf{y}\in A$ for which  $\mathbf{x}\prec \mathbf{y}$ we have $f(\mathbf{x})\ge f(\mathbf{y})$. We refer the reader to~\cite{Marshall_Olkin_Arnold} for further details and references on the theory of majorization.  

The \emph{Shannon entropy} of 
a random variable $X$, taking values in a finite set $\mathcal{S}_X$, is defined as  
\[
\bh(X) \,=\, -\sum_{x\in\mathcal{S}_X} \Prb(X=x)\cdot \log(\Prb(X=x)) \, ,
\]
where $\log$ denotes the binary logarithm, and with the usual convention that $0\log(0)=0$. We refer the reader to~\cite{Cover_Thomas} for a standard reference on the topic. 
The \emph{binary entropy function} is defined as $h(x)=-x\log(x)-(1-x)\log(1-x)$, for $x\in [0,1]$. Note that $h(x)$ is equal to the Shannon entropy of a Bernoulli random variable of mean $x\in [0,1]$, denoted $\text{Ber}(x)$. 
It is a well known fact that the binary entropy function is concave. 
Moreover, it is a well known fact (see~\cite[p.~101]{Marshall_Olkin_Arnold}) that 
the Shannon entropy of $X$ is  Schur-concave on the set $\{(p_1,\ldots,p_n)\in [0,1]^n : \sum_ip_i=1\}$, for all $n$. In other words, if $X$ and $Y$  are two random variables such that $p_X\prec p_Y$, then it holds 
$\bh(X)\ge \bh(Y)$. 

\subsection{Related work and main results}

In this article we shall be concerned with an entropic analogue of the Manickam-Mikl\'os-Singhi (MMS) conjecture, which is a statement about the sum in a sampling without replacement from a finite, zero-mean, population.  

More concretely, fix positive  integers $n\ge 2$ and $k<n$. Throughout the text, we denote by $\mathcal{N}_n$ the family    consisting of all 
multisets $P=\{x_1,\ldots,x_n\}$  of  real numbers whose sum is equal to zero. 
We view the elements of $\mathcal{N}_n$ as populations from which we sample $k$ elements without replacement. Formally, let  $P\in\mathcal{N}_n$ be fixed and choose   $\mbI\in \binom{[n]}{k}$  uniformly at random.  
Let   
\[
X_P = \sum_{i\in\mbI} x_i 
\]
be the sum of the $k$ elements that are sampled without replacement from $P$.  
What is a sharp lower bound on the probability that $X_P$ is non-negative? 

The following  has been conjectured almost forty years ago by Manickam and Mikl\'os~\cite{Manickam_Miklos} and by Manickam and Singhi~\cite{Manickam_Singhi}, driven by questions concerning the so-called first distribution invariant of the Johnson scheme. 

\begin{conjecture}[MMS, $1987$]\label{mms_conj}
If $n\ge 4k$, then for any $P\in\mathcal{N}_n$ it holds that  
\[
\Prb(X_P \ge 0) \,\ge\, \frac{k}{n} \, . 
\]
\end{conjecture}

If the conjecture is true, then the result is sharp as can be seen by considering the population $P_1 =\{1,\frac{-1}{n-1},\ldots,\frac{-1}{n-1}\}$, consisting of one element that is  equal to $1$ and $n-1$ elements 
that are all equal to $\frac{-1}{n-1}$, which clearly satisfies 
$\Prb(X_{P_1}\ge 0) = k/n$. 

One particular line of research, which has attracted considerable attention (see~\cite{AHS, Frankl, Pokrovskiy, Tyomkyn}, and references therein), focuses  (for fixed $k$) on upper bounds   on  the minimal integer $N_k$ having the property that for all $n\ge N_k$ and all  $P\in\mathcal{N}_n$ it holds $\Prb(X_P \ge 0) \ge k/n$. 
The MMS conjecture asserts that $N_k\le 4k$, and  
the best current upper bound on $N_k$,   due to Pokrovskiy~\cite{Pokrovskiy},  is $N_k \le 10^{46}k$. 
In other words, if $n\ge 10^{46}k$ then the MMS conjecture holds true. 

Let us remark that 
Conjecture~\ref{mms_conj} is not true without the assumption $n\ge 4k$. 
For instance, when $n=3k+1$,  the population 
$P_{3k-2}=\{\frac{1}{3k-2},\ldots,\frac{1}{3k-2},\frac{-1}{3}, \frac{-1}{3},\frac{-1}{3}\}$, consisting of $3k-2$ values that are equal to $\frac{1}{3k-2}$ and three values that are equal to $\frac{-1}{3}$, is such that it holds   $\Prb(X_{P_1}\ge 0) > \Prb(X_{P_{3k-2}}\ge 0)$, as can be easily verified, and so  the above-mentioned population $P_1$ is not minimal in this case. 
The following statement is a refinement of Conjecture~\ref{mms_conj} that covers all  non-trivial instances of the parameters $n$ and $k$. 

\begin{conjecture}[Aydinian, Blinovsky~\cite{Aydinian_Blinovsky}, $2012$]
\label{AB_conj}
Let $n\ge 2$ be an integer. Then 
for any $k\in [n-1]$ and any $P\in\mathcal{N}_n$ it holds  that 
\[
\Prb(X_P \ge 0) \,\ge\, \min_{i\in [n-1]} \Prb(H_i \ge ik/n) \, ,  
\]
where $H_i\sim\text{Hyp}(n,i,k)$. 
\end{conjecture}
The hypergeometric random variable $H_i$ in Conjecture~\ref{AB_conj} corresponds to the sum of a sampling  without replacement of size $k$ from the population 
\begin{equation}\label{special_P}
P_i = \left\{\frac{1}{i}, \ldots,\frac{1}{i}, \frac{-1}{n-i},\ldots,\frac{-1}{n-i}\right\} \, ,  
\end{equation}
consisting of $i\in [n-1]$ values that are equal to $\frac{1}{i}$ and $n-i$ values that are equal to $\frac{-1}{n-i}$. 
In this case, if we denote by $H_i\sim\text{Hyp}(n,i,k)$ the number of nonnegative elements in a sample of size $k$ without replacement from $P_i$, we have 
\[
X_{P_i} = \frac{1}{i}\cdot H_i - \frac{1}{n-i}\cdot (k-H_i) = \frac{n}{i(n-i)} \cdot \left( H_i - \frac{ik}{n}\right) \, ,
\]
and so $\Prb(X_{P_i}\ge 0) = \Prb(H_i \ge ik/n)$. 
In other words, provided the Aydinian-Blinovsky conjecture holds true, the populations $P_i, i\in [n-1]$, defined in~\eqref{special_P},  appear as extreme instances of the MMS conjecture. 
Note that Conjecture~\ref{AB_conj}, if true, reduces the problem of obtaining a sharp lower bound on $\Prb(X_P\ge 0)$ to the one of minimizing  the  ``tails"  $\Prb(H_i\ge ik/n)$ over $i\in [n-1]$.  Let us note that the latter optimization  problem can be solved under the assumption $n\ge 8k$ (see~\cite{Ai_Pelekis}) but, to the best of our knowledge, its solution is not known for smaller values of $n$.

In this article we investigate the above-mentioned question from an entropic point of view. 
Fix positive integers $n\ge 2$ and $k<n$, and let $P=\{x_1,\ldots,x_n\}\in\mathcal{N}_n$. 
If $X_P = \sum_{i\in \mathbb{I}}x_i$, where $\mathbb{I}\in\binom{[n]}{k}$ is chosen uniformly at random, 
what is a sharp lower bound on the Shannon entropy of $X_P$? 

To make the problem interesting, we have to assume that the population 
$P\in\mathcal{N}_n$ is non-trivial, i.e., not all elements $x_i$ are equal to zero. Under this assumption, we deduce the following entropic analogue of the MMS conjecture.

\begin{theorem}\label{main_thm}
Let $n\ge 2$ be an integer. Then 
for any $k\in [n-1]$ and any $P=\{x_1,\ldots,x_n\}\in\mathcal{N}_n$ such that $\sum_{i\in [n]}|x_i|>0$ it holds  that 
\[
\bh(X_P) \,\ge\, \bh(\text{Ber}(k/n)) \, ,
\]
where $\text{Ber}(k/n)$ is a Bernoulli random variable of mean $k/n$. 
The inequality becomes an equality when $P=\{1,\frac{-1}{n-1},\ldots,\frac{-1}{n-1}\}$.
\end{theorem}

The proof of Theorem~\ref{main_thm} relies on the following entropic analogue of the Aydinian-Blinovsky conjecture. 

\begin{theorem}\label{main_cor}
Let $n\ge 2$ be an integer. Then 
for any $k\in [n-1]$ and any $P=\{x_1,\ldots,x_n\}\in\mathcal{N}_n$ such that $\sum_{i\in [n]}|x_i|>0$ it holds  that 
\[
\bh(X_P) \,\ge\, \min_{i\in [n-1]} \, \bh(H_i) \, ,  
\]
where $H_i\sim\text{Hyp}(n,i,k)$. 
\end{theorem}

In other words, the populations $P_i, i\in [n-1]$, defined in~\eqref{special_P},  are extreme instances of the ``entropified MMS conjecture".
We provide two proofs of Theorem~\ref{main_cor}. 
The first proof  combines ideas from Sperner theory and the theory of majorization, and provides more than desired. 
Let us remark that an interplay between   Sperner theory  and the theory of majorization has been reported in~\cite{Madiman_Wang_Woo}, where  entropy inequalities are derived for sums of independent integer-valued
random variables. Our setting is different, as we are dealing with  sums of dependent random variables. 

A bit more precisely, our  approach is to consider the family  $\binom{[n]}{k}$ as a ranked   poset, which we refer to as the \emph{sign-split poset}, having the Strong Sperner Property. 
It turns out (see Lemma~\ref{obs1} below) that each value attained by the random variable $X_P$ corresponds to an antichain in this sign-split poset, whose Strong Sperner  Property (Lemma~\ref{obs2})  readily implies    that the probability mass function, $p_{X_P}$, of $X_P$ is majorized by the probability mass function, $p_{H_i}$, of a hypergeometric random variable $H_i\sim\text{Hyp}(n,i,k)$, for some $i\in [n-1]$ (Lemma~\ref{maj_lemma}). 
This guarantees that, for any Schur-concave function $f$, it holds 
$f(p_{X_P})\ge f(p_{H_i})$. In particular, the latter holds true for the  Shannon entropy and provides the first proof of Theorem~\ref{main_cor}.   
As a byproduct of the aforementioned   majorization between the   probability mass functions, we deduce  the  following Littlewood-Offord-type result for  $X_P$.  

\begin{theorem}\label{main_thm2}
For any $k<n$ and any $P=\{x_1,\ldots,x_n\}\in\mathcal{N}_n$ such that $\sum_{i\in [n]}|x_i|>0$ it holds that 
\[
\sup_{x\in \mathbb{R}}\,\Prb(X_P = x) \,\le\, \max_{i\in [n-1]} \,\Prb(H_i = m_i) \, ,
\]
where $H_i\sim\text{Hyp}(n,i,k)$ and $m_i$ is a mode of $H_i$. 
\end{theorem}

We also provide a second proof of Theorem~\ref{main_cor}, which is of a more elementary nature since it   avoids completely Sperner theory and Schur-concavity. It exploits the appealing  properties of the Shannon entropy, and draws inspiration from  Lubell's well-known proof of the LYM inequality; so it may be of independent interest.

Finally, to prove Theorem~\ref{main_thm}, observe that  Theorem~\ref{main_cor}  reduces the problem of minimizing $\bh(X_P)$ to the one of minimizing the entropies $\bh(H_i)$ over $i\in [n-1]$ and, unlike the setting of the MMS conjecture, we are able to solve the latter optimization problem  (Lemma~\ref{last_lemma}) 
using a coupling argument and the appealing properties of the hypergeometric law. 

\subsection{Organisation} 

The remaining part of our article is organised as follows. 
In Section~\ref{sec:sperner} we collect the basic definitions and results from Sperner theory that will be need in the proofs of our main results. 
In Section~\ref{sec:sign_split} we define the sign-split poset of a population in $\mathcal{N}_n$ and prove that it possesses the Strong Sperner Property.  Our main results are proven in   
Section~\ref{proof_of_main_thm}, and in Section~\ref{Lubell_proof} we provide a second proof of Theorem~\ref{main_cor}.  
Our article ends with Section~\ref{concluding_remarks}  in which we discuss how our approach may be employed in similar  settings.

\section{The sign-split poset}

\subsection{Preliminaries from Sperner theory}
\label{sec:sperner}

As already mentioned, 
our main results employ ideas from Sperner theory of partial ordered sets (posets), which we briefly discuss here and refer the reader to~\cite{Engel} for an excellent textbook on the topic.

Recall that a \emph{finite poset} is a pair $Q=(X,\le)$, where $X$ is a finite set equipped with a partial order relation $\le$ which is  reflexive, antisymmetric and transitive. We write $p\le q$ for two comparable elements $p,q$ of the poset, and $p < q$ when $p\le q$ and $p\neq q$. 
An element $q\in X$ is said to \emph{cover} the elements $p\in X$, and is denoted $p\lessdot q$, if 
$p<q$ and $p < q'\le q$ implies $q=q'$. 
An element $p\in X$ is called \emph{maximal} if $p\le q$ implies $q=p$. 
Similarly, an element $p\in X$ is called \emph{minimal} if $q\le p$ implies $q=p$.

A \emph{chain} in a poset $Q$ is a subset of pairwise comparable elements, i.e., a set $C\subseteq X$ of the form 
$C = (p_1 < p_2< \cdots < p_m)$, for some integer $m\ge 1$. 
The integer $m$ is the \emph{length} of the chain. 
A chain $C$ in a poset $P$ is called \emph{saturated} if it has the form $C = (p_1 \lessdot p_2\lessdot \cdots \lessdot p_m)$. 
  
Given a positive integer $\ell\ge 1$, an $\ell$-\emph{Sperner family} in a poset $Q$ is a subset $S\subseteq X$ such that $|S\cap C|\le \ell$, for all chains $C$ of $Q$.  A $1$-Sperner family is referred to as an \emph{antichain}. 
In other words, an $\ell$-\emph{Sperner family} is a subset of $X$ which does not contain a chain of length $\ell+1$. 

A chain $C$ in a poset $Q=(X,\le)$ is \emph{maximal} if there is no chain $C'$ in $Q$, which is different from $C$,  such that $C\subseteq C'$. 
If every maximal chain of the poset has the same length, then the poset is called \emph{graded}. In this case we can define a unique \emph{rank function}  $r:X\to \{0,1,\ldots\}$ from $X$ to the natural numbers which satisfies   
$r(q)=r(p)+1$ whenever $p\lessdot q$, and is such that  every minimal element has \emph{the same} rank and every maximal element has \emph{the same} rank. The \emph{rank} of an element $p\in X$ is the value $r(p)$, and the \emph{rank of the poset} is defined as $r(Q)= \max\{r(p): p\in X\}$. The \emph{minimum rank} of $Q$ is the integer 
$\mu(Q) = \min\{r(p): p\in X\}$. A \emph{ranked} poset is a poset in which a rank function can be defined. 

The $i$-\emph{th level set} of a ranked poset $Q=(X,\le)$ is the set 
$N_i(Q) = \{p\in X : r(p)=i\}$, for integer $i\ge 0$. 
Note that, by definition, each element $p\in X$ belongs to the level set $N_{r(p)}(Q)$. 
The $i$-\emph{th Whitney number} is the size of the $i$-th level set, denoted $W_i(Q)$. Hence $W_i(Q)= |N_i(Q)|$. 

A ranked poset $Q=(X,\le)$ is said to have the 
$\ell$-\emph{Sperner property} if the maximum size of a 
$\ell$-Sperner family in $Q$ equals the sum  of the \emph{$\ell$ largest Whitney numbers} of $Q$. The ranked poset $Q$ has the \emph{Strong Sperner property} if it has the $\ell$-Sperner Property for all integers  $\ell\ge 1$. 

Given an element $p$ of a poset $Q=(X,\le)$, we define the \emph{upper shadow} of $p$ to be the set 
\[
\nabla(p) = \{q\in X : p\lessdot q\} \, . 
\]
The \emph{upper shadow} of a set $A\subseteq X$ is defined  as  $\nabla(A)=\{q\in X : p\lessdot q, \text{ for some } p\in A\}$.

A ranked poset $Q=(X,\le)$, with rank function $r$, is said to have the 
\emph{normalized matching property (NMP)} if it satisfies 
\[
\frac{|A|}{W_i(Q)}\le \frac{|\nabla (A)|}{W_{i+1}(Q)} \, ,\, \text{ for all } A\subseteq N_i \, \text{ and all } \, i \in \{\mu(Q),\ldots, r(Q)-1 \}\, .
\]
Briefly, a ranked poset with the NMP property is referred to as a \emph{normal poset}. The following well-known result will be employed in the proof of Theorem~\ref{main_cor}. 

\begin{theorem}\label{engel_thm}
Let $Q=(X,\le)$ be a normal poset. Then $Q$ has the Strong Sperner Property. 
\end{theorem}
\begin{proof}
See~\cite[Theorem~4.5.1 and Corollary~4.5.3]{Engel}.  
\end{proof}

\subsection{Normality of  the sign-split  poset}\label{sec:sign_split}

In this section we associate a ranked poset to each non-trivial  population  from 
$\mathcal{N}_n$. More precisely, let $n\ge 2$ and $k\in [n-1]$ be fixed integers, and consider a multiset  $P=\{x_1,\ldots,x_n\}$  such that $\sum_{i\in [n]}x_i=0$ and $\sum_{i\in [n]}|x_i| >0$. 
Consider the following partition of the index set:
\[
A = \{i\in [n] : x_i \ge 0\} \quad\text{ and }\quad B = \{i\in [n]: x_i<0\} \, .
\]
Note that every element $F\in\binom{[n]}{k}$ can be partitioned  as 
$F = F_{A} \cup F_B$, where $F_A = F\cap A$ and $F_B = F\cap B$. 
We define a partial order on the collection  $\binom{[n]}{k}$ via the relation 
\[
F\le G \quad \iff \quad F_A\subseteq G_A \,\text{ and }\, G_B\subseteq F_B \, .
\]
Then $Q=(\binom{[n]}{k}, \le)$ is a poset, which will be referred to as the \emph{sign-split  poset} corresponding to  $P$.  
We turn $Q$ into a ranked poset by defining the rank of an element $F\in \binom{[n]}{k}$ to be equal to 
\[
r(F) = |F_A| \, .
\]
The $i$-th level of $Q$ is the set $N_i:=\{F\in\binom{[n]}{k} : |F_A| = i\}$,   and therefore the   $i$-th Whitney number of the sign-split poset is equal to 
\[
W_i := \binom{|A|}{i}\binom{n-|A|}{k-i} \, , \, \text{ for } \, i\in\big\{\max\{0,k-(n-|A|)\}, \ldots, \min\{k,|A|\}-1 \,\big\}\, .
\]
The proofs of our main results are based on the following observations. 

\begin{lemma}\label{obs1}
Let $n\ge 2$ be an integer and fix $k\in [n-1]$. 
Let also $P=\{x_1,\ldots,x_n\}$ be such that  $\sum_{i\in [n]}x_i=0$ and $\sum_{i\in [n]}|x_i| >0$. 
Then, for any real number $x$, the collection  
\[
\mathcal{F}_x = \left\{F\in\binom{[n]}{k} \,:\, \sum_{i\in F} x_i = x \right\}
\]
is an antichain of the sign-split  poset corresponding to $P$. 
\end{lemma}
\begin{proof}
Suppose that the conclusion is not true. Then there exist $F,G\in\mathcal{F}_x$ such that $F<G$. 
Hence, using the fact that $x_i<0$, for $i\in B$, we deduce that 
\[
0 = \sum_{i\in G} x_i - \sum_{i\in F}x_i = \sum_{i\in G_A\setminus F_A} - \sum_{i\in G_B\setminus F_B} x_i  \,>\, 0  \, ,
\]
which is a contradiction. 
\end{proof}

The following result is the main ingredient of our approach. 

\begin{lemma}\label{obs2}
The sign-split poset $Q=(\binom{[n]}{k}, \le)$ has the Strong Sperner Property. 
\end{lemma}
\begin{proof}
It follows from Theorem~\ref{engel_thm} that it is enough to show that $Q$ is normal. To this end, we have to show that 
for all $i \in \{\max\{0,k-(n-|A|)\}, \ldots, \min\{k,|A|\}-1\}$ and all $\mathcal{F}\subseteq N_i$ it holds 
\begin{equation}\label{double_count}
\frac{|\mathcal{F}|}{W_i}\le \frac{|\nabla(\mathcal{F})|}{W_{i+1}} \, .
\end{equation}
Let $i$ and $\mathcal{F}$ be as above. We count pairs 
$(F,G)$ such that $F\in \mathcal{F}$, $G\in\nabla(\mathcal{F})$ and $F\lessdot G$. 
For every $F\in \mathcal{F}$ there are exactly $(|A|-i)\cdot (k-i)$ sets $G\in\nabla(\mathcal{A})$ such that $F\lessdot G$; those obtained by replacing one element of $F_B$ with one element from $A\setminus F_A$. 
Hence the total number of pairs $(F,G)$ is equal to 
$|\mathcal{F}|\cdot (|A|-i)\cdot (k-i)$. 

On the other hand, for every $G\in\nabla(\mathcal{F})$ there are at most $(i+1)\cdot (n-|A|-k+i+1)$ elements $F\in\mathcal{F}$ such that $F\lessdot G$; those obtained by replacing an element of $G_A$ with an element from $B\setminus G_B$. 
Hence it holds 
\[
|\mathcal{F}|\cdot (|A|-i)\cdot (k-i) \le |\nabla(\mathcal{F})|\cdot (i+1)\cdot (n-|A|-k+i+1) \, .
\]
The latter inequality  and  the fact that $W_i=\binom{|A|}{i}\binom{n-|A|}{k-i}$ imply that the desired inequality~\eqref{double_count} holds true, and  the result follows. 
\end{proof}

\section{Proofs of the main results}
\label{proof_of_main_thm}

Throughout this section, we let the integers $n\ge 2$ and $k\in [n-1]$ be fixed, and we also fix a population 
$P=\{x_1,\ldots,x_n\}\in\mathcal{N}_n$ which satisfies  
$\sum_{i\in [n]}|x_i|>0$. 

Choose $\mathbb{I}\in\binom{[n]}{k}$ uniformly at random and set $X_P = \sum_{i\in\mathbb{I}}x_i$. 
Let $\mathcal{X}_P=\{s_1,\ldots,s_m\}$ be the set of real numbers for which it holds $\Prb(X_P=s_j)>0$, for $j\in [m]$.  Note that $m\ge 2$. 
For $j\in [m]$, let 
\[
\mathcal{F}_j = \left\{F\in \binom{[n]}{k} \,:\, \sum_{l\in F}x_l = s_j \right\} \quad \text{ and }\quad f_j = |\mathcal{F}_j| \, . 
\]
Then it holds 
\[
\Prb(X_P = s_j) = \frac{f_j}{\binom{n}{k}} \, , \,\text{ for } \, j \in [m] \, .
\]
Let $p_{X_P} = (f_1/\binom{n}{k},\ldots, f_m/\binom{n}{k})$ be the probability mass function of $X_P$.  
For $i\in [n-1]$, let $H_i\sim\text{Hyp}(n,i,k)$ and 
let $p_{H_i}$ be its probability mass function.  
In the following lemma we show that the probability mass function of $X_P$ is majorized by the probability mass function of a hypergeometric random variable. 

\begin{lemma}\label{maj_lemma}
There exists $i\in [n-1]$ such that $p_{X_P} \prec p_{H_i}$.  
\end{lemma}
\begin{proof}
Let $Q=(\binom{[n]}{k}, \le)$ be the sign-split poset corresponding to $P$. 
Recall that in the definition of the sign-split we partition the indices of the elements in the population into the sets
\[
A=\{j\in [n]: x_j \ge 0\} \quad\text{ and }\quad B = \{j\in [n] : x_j < 0\} \, .
\]
Set $i=|A|$, and note that $i\in [n-1]$. 
Let
$a=\max\{0,k-(n-i)\}$ and $b=\min\{k,i\}$, and observe that the non-empty levels of $Q$ are $N_a,\ldots,N_b$, where $N_j = \{F\in \binom{[n]}{k} : |F\cap A| =j \}$. 
Then the Whitney numbers of $Q$  
are equal to 
\[
W_j = |N_j| = \binom{i}{j}\binom{n-i}{k-j} \, , \, \text{ for } \, j\in \{a,\ldots,b\}\, .
\]
Set $d = b-a+1$, and let $W_{(1)}\ge \cdots\ge W_{(d)}$ be the Whitney numbers of $Q$ in decreasing order.  
Since every maximal saturated chain of
$Q$ has one element in each level, and the sum $\sum_{r\in F}x_r$  strictly increases along chains, it follows that 
the random variable $X_P$ takes at least $d$ distinct values. Thus $m\ge d$.  We extend the
sequence $W_{(1)},\ldots,W_{(d)}$ by setting \(W_{(l)}=0\) for $d<l\leq m$.  Let $\mathbf{w}_i$ denote the vector $(W_{(1)},\ldots,W_{(m)})$.

Now let $f_{(1)}\ge \cdots \ge f_{(m)}$ be the coordinates of the  vector
$\mathbf{v}_P = (f_1,\ldots,f_m)$ in decreasing order.  
Note that it holds $\sum_{j\in [m]}f_{(j)}=\sum_{j\in [m]}W_{(j)}=\binom{n}{k}$.
We now claim that,  for all 
$\ell\in [m-1]$, it holds 
\[
\sum_{j=1}^{\ell} f_{(j)} \le \sum_{j=1}^{\ell} W_{(j)} \, .
\]
To see the claim, note that  Lemma~\ref{obs1} implies that each $\mathcal{F}_j$ is an antichain of the sign-split poset. Fix $\ell\in [m-1]$, 
and assume that the integers $f_{(1)}\ge\ldots \ge f_{(\ell)}$ are the cardinalities of the families $\mathcal{F}_{s_1},\ldots,\mathcal{F}_{s_{\ell}}\subseteq \binom{[n]}{k}$, respectively. 
Then the collection $\mathcal{G} = \mathcal{F}_{s_1}\cup \cdots\cup \mathcal{F}_{s_{\ell}}$ is a subset of the sign-split poset which does not contain a chain of size $\ell+1$.  
Since the sign-split poset has the Strong Sperner property, by Lemma~\ref{obs2}, the claim follows. 

In particular, we conclude that the vector 
$\frac{1}{\binom{n}{k}}\cdot \mathbf{v}_P$ is majorized by the vector $\frac{1}{\binom{n}{k}}\cdot \mathbf{w}_i$. 
The result follows upon observing the the former vector is the probability mass function of $X_P$, while the latter is the probability mass function of $H_i\sim\text{Hyp}(n,i,k)$. 
\end{proof}

Note that Lemma~\ref{maj_lemma} immediately implies Theorem~\ref{main_thm2}. Moreover, since the Shannon entropy is a Schur-concave function of the probability mass function, Theorem~\ref{main_cor}  follows immediately  from  Lemma~\ref{maj_lemma}, 
and we are thus left with proving Theorem~\ref{main_thm}. 
The following lemma, combined with Theorem~\ref{main_cor},  completes the proof of  Theorem~\ref{main_thm}. 

\begin{lemma}\label{last_lemma}
Let $n\ge 2$ and $k\in [n-1]$. Then for all $i\in [n-1]$ it holds 
\[
\bh(H_i) \ge \bh(\text{Ber}(k/n)) \, ,
\]
where $H_i\sim\text{Hyp}(n,i,k)$ and 
$\text{Ber}(k/n)$ is a Bernoulli random variable of mean $k/n$. 
\end{lemma}
\begin{proof} 
The result is clearly true when $n=2$, so we may assume that $n\ge 3$. 
Note that the probability mass function of $H_1$ is equal to $(\frac{k}{n}, 1 - \frac{k}{n})$, and therefore it holds 
$\bh(H_1)=\bh(\text{Ber}(k/n))$. 
Note also that  for each $i\in [n-1]$ it holds 
$k-H_i \sim\text{Hyp}(n,n-i,k)$. Hence $k-H_i \sim H_{n-i}$ and so we have $\bh(H_i) = \bh(H_{n-i})$. 
In particular, we may assume that $i\le \lfloor n/2\rfloor$, and it  is therefore  enough to show that 
\begin{equation}\label{enough1}
\bh(H_{i+1}) \ge \bh(H_i) \, , \, \text{ for } \, i \le \lfloor n/2\rfloor -1 \, .
\end{equation}
We prove~\eqref{enough1} via a coupling argument.  
More concretely, let $\mathbb{I}\in\binom{[n]}{k}$ be chosen uniformly at random and, for each $i\in [n-1]$, consider the random variable  
\[
Z_i = |\mathbb{I} \cap [i]| \, .
\]   
Note that $H_i \sim Z_i$, for all $i\in [n-1]$, and therefore it is enough to show that 
\begin{equation}\label{enough1.2}
\bh(Z_{i+1}) \ge \bh(Z_i) \, , \, \text{ for } \, i \le \lfloor n/2\rfloor -1 \, .
\end{equation}
Observe that it holds 
\[
Z_{i+1} = Z_i + \mathbf{1}_{\{i+1\in \mathbb{I}\}} \, ,
\]
where $\mathbf{1}_{\{i+1\in \mathbb{I}\}}$ denotes the indicator of the event $\{i+1\in \mathbb{I}\}$. 
In other words, if we know the pair $(Z_i, \mathbf{1}_{\{i+1\in \mathbb{I}\}})$ then we can recover $Z_{i+1}$, and if we know the pair 
$(Z_{i+1}, \mathbf{1}_{\{i+1\in \mathbb{I}\}})$ then we can recover $Z_i$. 
This implies that 
\[
\bh(Z_i, \mathbf{1}_{\{i+1\in \mathbb{I}\}}) = \bh(Z_{i+1}, \mathbf{1}_{\{i+1\in \mathbb{I}\}}) 
\]
which, using the chain rule, can be equivalently written as 
\[
\bh(Z_i) + \bh(\mathbf{1}_{\{i+1\in \mathbb{I}\}}\mid Z_i) = 
\bh(Z_{i+1}) + \bh(\mathbf{1}_{\{i+1\in \mathbb{I}\}}\mid Z_{i+1}) \, .
\]
Hence, in order to prove~\eqref{enough1.2}, it is enough to show that 
\begin{equation}\label{enough2}
\bh(\mathbf{1}_{\{i+1\in \mathbb{I}\}}\mid Z_i) \ge \bh(\mathbf{1}_{\{i+1\in \mathbb{I}\}}\mid Z_{i+1}) \, , \, \text{ for } \, i \le \lfloor n/2\rfloor -1 \, .
\end{equation}
Let $i\le \lfloor n/2\rfloor -1$ be a positive integer. Note that 
\[
\Prb(\mathbf{1}_{\{i+1\in\mathbb{I}\}}=1 \mid Z_i = j)= \frac{\binom{i}{j}\binom{n-i-1}{k-j-1}/\binom{n}{k}}{\binom{i}{j}\binom{n-i}{k-j}/\binom{n}{k} }
= \frac{k-j}{n-i} \, .
\]
Hence, if  $h(\cdot)$ denotes the binary entropy function, we deduce that
\begin{eqnarray*}
\bh(\mathbf{1}_{\{i+1\in \mathbb{I}\}}\mid Z_i) &=& \sum_{j} \bh(\mathbf{1}_{\{i+1\in \mathbb{I}\}}\mid Z_i=j) \cdot \Prb(H_i=j)\\
&=& \sum_{j} h\left( \frac{k-j}{n-i}\right) \cdot \Prb(Z_i=j) \\
&=& \Exp\left(h\left(\frac{k-Z_i}{n-i} \right) \right) \\
&=& \Exp\left(h\left(\frac{Z_{n-i}}{n-i} \right) \right) \, ,
\end{eqnarray*}
where the last inequality follows from the fact that $k-Z_i\sim k-H_i\sim H_{n-i}\sim Z_{n-i}$. 
Similarly, we have 
\[
\Prb(\mathbf{1}_{\{i+1\in\mathbb{I}\}}=1 \mid Z_{i+1} = \ell)= \frac{\binom{i}{\ell-1}\binom{n-i-1}{k-\ell}/\binom{n}{k}}{\binom{i+1}{\ell}\binom{n-i-1}{k-\ell}/\binom{n}{k} }
= \frac{\ell}{i+1} \, , 
\]
and therefore we may write 
\[
\bh(\mathbf{1}_{\{i+1\in \mathbb{I}\}}\mid Z_{i+1}) = \Exp\left(h\left(\frac{Z_{i+1}}{i+1}\right)\right) \, .
\]
Summarizing the above, it follows from~\eqref{enough2} that it is enough to show that 
\begin{equation}\label{enough3}
\Exp\left(h\left(\frac{Z_{n-i}}{n-i} \right) \right) \ge \Exp\left(h\left(\frac{Z_{i+1}}{i+1}\right)\right) \, , \, \text{ for } \, i \le \lfloor n/2\rfloor -1 \, .
\end{equation}
We prove a bit more. Namely, we show that the sequence 
\[
E_m := \Exp\left(h\left(\frac{Z_{m}}{m}\right)\right) \, , \text{ for } \, m=1,\ldots,n-1
\]
is non-decreasing. The desired inequality~\eqref{enough3} then follows upon observing that the assumption $i \le \lfloor n/2\rfloor -1$ implies that $n-i\ge i+1$. To this end, let $m\in [n-2]$ and note that, since $Z_{m+1} = Z_m + \mathbf{1}_{\{m+1\in \mathbb{I}\}}$, we may write 
\begin{eqnarray*}
\Exp\left(\frac{Z_m}{m} \mid Z_{m+1}=\ell \right) &=&  
\Exp\left(\frac{Z_{m+1}-\mathbf{1}_{\{m+1\in\mathbb{I}\}}}{m} \mid Z_{m+1}=\ell \right)\\
&=& \frac{\ell}{m} - \frac{1}{m}\cdot \Prb(\mathbf{1}_{\{m+1\in\mathbb{I}\}}= 1 \mid Z_{m+1}=\ell)\\
&=& \frac{\ell}{m} -\frac{1}{m}\cdot \frac{\ell}{m+1}\\
&=& \frac{\ell}{m+1} \, .
\end{eqnarray*}
Since the binary entropy function is concave, the conditional  Jensen inequality yields 
\[
\Exp\left(h\left(\frac{Z_m}{m}\right) \mid Z_{m+1}=\ell \right) \,\le\, h\left(\Exp\left(\frac{Z_m}{m} \mid Z_{m+1}=\ell\right)\right) = h\left( \frac{\ell}{m+1}\right) \, .
\]
In particular, this implies that 
\[
E_m=\Exp\left(h\left(\frac{Z_m}{m}\right)\right) = 
\Exp\left( \Exp\left(h\left(\frac{Z_m}{m}\right) \mid Z_{m+1}  \right)\right) \le \Exp\left(h\left(\frac{Z_{m+1}}{m+1}\right)  \right)= E_{m+1} \, ,
\]
as desired.  
\end{proof}

\section{An elementary proof of Theorem~\ref{main_cor}}
\label{Lubell_proof}

In this section we provide a proof of Theorem~\ref{main_cor} that avoids Sperner theory and majorization.  The proof draws inspiration from Lubell's proof for the LYM inequality.
Briefly, the idea is to show that drawing from the  random variable $X_P$ 
is equivalent to first choosing  a  maximal chain of the sign-split poset uniformly at random, then choosing a $k$-set of the selected chain with probability proportional to its rank and finally  adding the elements 
of this $k$-set. 
It turns out that the function that maps the rank of the selected $k$-set  to the sum of its elements is  injective, and  allows to bypass the use of Schur-convexity. 
The Sperner theory  is bypassed upon observing that a uniformly chosen maximal chain in the sign-split poset is completely determined by two random permutations of the sign-split sets.

\begin{proof}[Second proof of Theorem~\ref{main_cor}]
Consider a partition the indices of the elements in the population into non-negative and negative elements:
\[
A=\{j\in [n]: x_j \ge 0\} \quad\text{ and }\quad B = \{j\in [n] : x_j < 0\} \, .
\]
Set $i=|A|$, and note that $i\in [n-1]$. 
Now choose a random permutation of the elements in $A$, say $a_1,\ldots,a_i$, and, independently, a random permutation of the elements in $B$, say $b_1,\ldots,b_{n-i}$. 
For every $j\in\{\max\{0,k-(n-i)\},\dots,\min\{k,i\}\}$, consider the (random) set $C_j = A_j \cup B_{k-j}$, where 
\[
A_j = \{a_1,\ldots,a_j\} \quad \text{ and } \quad B_{k-j}=\{b_1,\ldots,b_{k-j}\} \, .
\]
Here, we set $A_0=B_{0}=\emptyset$.  
Now consider the (random) collection $\mathcal{C}= \{C_j\}_{j=\max\{0,k-(n-i)\}}^{\min\{k,i\}}$. 
Observe that each set $C_j$ is uniformly distributed on the collection 
$\mathcal{B}_j$ consisting of all $F\in\binom{[n]}{k}$ for which 
$|F\cap A|=j$. 

We claim that the map $f: \{\max\{0,k-(n-i)\},\dots,\min\{k,i\}\}\to \mathbb{R}$ with 
$f(j)= \sum_{l\in C_j}x_l$ is strictly increasing. Indeed, since $C_{j+1} = (C_j\setminus \{b_{k-j}\})\cup \{a_{j+1}\}$, it holds that 
\[
f(j+1)-f(j)=\sum_{l\in C_{j+1}}x_l - \sum_{l\in C_j}x_l = x_{a_{j+1}} - x_{b_{k-j}} > 0 \, .
\]
Now let $H_i\sim\text{Hyp}(n,i,k)$ be independent of the collection $\mathcal{C}$, and note that  for a given $F\in\mathcal{B}_j$ it holds 
\begin{eqnarray*}
\Prb(C_{H_i} = F) &=& \Prb(C_{H_i}= F \mid H_i=j)\cdot \Prb(H_i=j)\\ 
&=&  \Prb(C_{j}= F)\cdot \Prb(H_i=j)
 = \frac{1}{\binom{i}{j}\binom{n-i}{k-j}}\cdot \frac{\binom{i}{j}\binom{n-i}{k-j}}{\binom{n}{k}} = \frac{1}{\binom{n}{k}} \, .
\end{eqnarray*}
This implies that $C_{H_i}$ is uniformly distributed on $\binom{[n]}{k}$, and therefore we conclude that  
\[
X_P \,\sim\, f(H_i)= \sum_{l \in C_{H_i}}x_l \, .
\]
Since $f$ is injective and $H_i$ is independent of $\mathcal{C}$, we deduce that 
\[
\bh(f(H_i)\mid \mathcal{C}) = \bh(H_i\mid \mathcal{C}) = \bh(H_i) \, .
\]
Hence, since conditioning does not increase entropy, it holds that  
\[
\bh(X_P) = \bh(f(H_i)) \ge \bh(f(H_i)\mid \mathcal{C}) = \bh(H_i) \, ,
\]
and the result follows. 
\end{proof}

\section{Concluding remarks}\label{concluding_remarks}

In this article we presented an ``entropified version" of the MMS conjecture. Our approach is based on Lemma~\ref{maj_lemma}, which guarantees that the probability mass function of $X_P$ is majorized by the probability mass function of a hypergeometric random variable $H_i\sim\text{Hyp}(n,i,k)$. The Schur-concavity of the Shannon entropy then  implies that an ``entropified version" of the Aydinian-Blinovsky holds true,  
which reduces the problem to the one of minimizing 
the Shannon entropy of the extremal random variables $H_i$, for $i\in [n-1]$. 

Let us remark that our approach may be employed for any Schur-concave function of the probability mass function. For instance, the \emph{R\'enyi entropy} of a random variable $X$, taking values in the set 
$\mathcal{S}_X$ and having probability mass function $p_X$, is defined as 
\[
\bh_{\alpha}(X) = \frac{1}{1-\alpha}\cdot\log\left(\sum_{x\in\mathcal{S}_X} (p_X(x))^{\alpha} \right) \, , \, \text{ for } \, \alpha \in (0,\infty) \, .
\]
It is a well knwon fact that the R\'enyi entropy is a Schur-concave function of the probability mass function, and therefore  Lemma~\ref{maj_lemma} yields the following analogue of the Aydinian-Blinovsky conjecture in this setting.

\begin{corollary}\label{Renyi_entropy}
Let $n\ge 2$ be an integer. Then 
for any $k\in [n-1]$ and any $P=\{x_1,\ldots,x_n\}\in\mathcal{N}_n$ such that $\sum_{i\in [n]}|x_i|>0$ it holds that 
\[
\bh_{\alpha}(X_P) \,\ge\, \min_{i\in [n-1]}\,\bh_{\alpha}(H_i) \, ,
\]
where $H_i\sim\text{Hyp}(n,i,k)$. 
\end{corollary}

We are unable to solve the optimization problem $\min_{i\in [n-1]}\,\bh_{\alpha}(H_i)$, but our numerical experiments suggest that $\bh_{\alpha}(H_1)$ is minimal.

\section*{Acknowledgements}

JA and OK were funded by the Czech Science Foundation project 24-11820S (``Automatic Combinatorialist").


\begin{thebibliography}{99}


\bibitem{Ai_Pelekis} J.~Ai, C.~Pelekis, \textit{On lower bounds for hypergeometric tails}, Statistics \& Probability Letters \textbf{236} (2026) Article 110798.  

\bibitem{AHS} N.~Alon, H.~Huang, B.~Sudakov, 
\textit{Nonnegative $k$-sums, fractional covers, and probability of small deviations}, 
J. Combin. Theory Ser. B \textbf{102} (2012) 784--796.

\bibitem{Aydinian_Blinovsky} H.~Aydinian, V.~Blinovsky, \textit{A Remark on the Problem of Nonnegative k-Sums}, Problems of Inform. Transmission \textbf{48}(4) (2012) 56--61. 

%\bibitem{Bier}  T. Bier, \textit{A distribution invariant for the association schemes and strongly regular graphs}, Linear
algebra and its applications \textbf{57} (1984) 230--2521.

%\bibitem{Bier_Manickam} T.~Bier, N.~Manickam, \textit{The first distribution invariant of the Johnson scheme}, SEAMS Bull. Math. \textbf{11} (1987) 61--68.


\bibitem{Cover_Thomas} T.M.~Cover, J.A.~Thomas, \textit{Elements of information theory}, 2nd Edition, John Wiley \& Sons, 2005. 

\bibitem{Engel} K.~Engel, \textit{Sperner Theory}, Cambridge University Press, 1997.

\bibitem{Frankl} P.~Frankl, \textit{On the number of nonnegative sums
}, J. Combin. Theory Ser. B
\textbf{103}(5) (2013) 647--649.

\bibitem{Madiman_Wang_Woo} M.~Madiman, L.~Wang, and J.O.~Woo, \textit{Majorization and Renyi entropy inequalities via Sperner
theory}, Discrete Mathematics \textbf{342}(10) (2019) 2911--2923.


\bibitem{Manickam_Miklos} N.~Manickam, D.~Mikl\'os, \textit{On the number of non-negative partial sums of a non-negative
sum}, Colloq. Math. Soc. Janos Bolyai \textbf{52} (1987) 385--392.

\bibitem{Manickam_Singhi} N.~Manickam, N.~M.~Singhi, \textit{First distribution invariants and EKR theorems}, J. Combin. Theory Ser. A \textbf{48} (1988) 91--103.

\bibitem{Marshall_Olkin_Arnold} A.~W.~Marshall, I.~Olkin, and B.~C.~Arnold, \textit{Inequalities: Theory of Majorization and Its
Applications}, 2nd ed., Springer Science+Business Media, LLC, New York, (2011).





\bibitem{Pokrovskiy} A.~Pokrovskiy, \textit{A linear bound on the Manickam–Mikl\'os–Singhi conjecture}, J. Combin. Theory Ser. A
\textbf{133} (2015) 280--306.


\bibitem{Tyomkyn}  M.~Tyomkyn, \textit{An improved bound for the Manickam-Mikl\'os-Singhi conjecture}, 
European Journal of Combinatorics \textbf{33}(1) (2012)
27--32. 







 

\end{thebibliography}
\end{document}